\documentclass{article}
\usepackage{amsmath}

\usepackage{amssymb}
\usepackage{latexsym}
\usepackage{latexsym,boxedminipage}
\usepackage{epsfig}

\newcommand{\la}{\lambda}

\newcommand{\qed}{\hfill{$\Box$}}

\newcommand{\A}{{\mathcal A}}
\newcommand{\B}{{\mathcal B}}
\newcommand{\C}{{\mathcal C}}
\newcommand{\D}{{\mathcal D}}
\newcommand{\E}{{\mathcal E}}
\newcommand{\F}{{\mathcal F}}
\newcommand{\G}{{\mathcal G}}
\newcommand{\mH}{{\mathcal H}}

\newcommand{\K}{{\mathcal K}}
\newcommand{\mL}{{\mathcal L}}

\newcommand{\mP}{{\mathcal P}}

\newtheorem{lemma}{Lemma}
\newtheorem{proposition}{Proposition}
\newtheorem{theorem}{Theorem}

\title{Five Guidelines for Partition Analysis with Applications to
Lecture Hall-type Theorems
}

\author{Sylvie Corteel\\
{CNRS PRiSM, UVSQ}\\
{45 Avenue des Etats-Unis}\\
{ 78035 Versailles, France}\\
{\tt syl@prism.uvsq.fr}\\
\and
Sunyoung Lee
\thanks{Research supported in part by NSF grant DMS-0300034}\\
{Computer Science}\\
{N. C. State University }\\
{ Raleigh, NC 27695}\\
{\tt slee7@unity.ncsu.edu}
\and
Carla D. Savage
\thanks{Research supported in part by NSF grants DMS-0300034 and INT-0230800}\\
{Computer Science}\\
{N. C. State University }\\
{ Raleigh, NC 27695}\\
{\tt savage@csc.ncsu.edu}
}

\date{August 23, 2005; revised May 19, 2006}

\begin{document}
\maketitle
\begin{abstract}
Five simple guidelines are proposed to compute the generating function
for the nonnegative integer solutions of a system of linear 
inequalities.  In contrast to other approaches, the emphasis is on
deriving recurrences.  We show how to use the guidelines strategically 
to solve some nontrivial enumeration problems in the theory of partitions and
compositions. This includes a strikingly different approach to lecture hall-type
theorems, with  new $q$-series identities arising in the process.
For completeness, we prove that the guidelines suffice to find the
generating function for 
any system of homogeneous linear inequalities with
integer coefficients.
  The guidelines can be viewed as a simplification of
MacMahon's partition analysis with ideas from
matrix
techiniques, Elliott reduction, and  ``adding a slice''.
\end{abstract}

\section{Introduction}

This continues our work in
\cite{ineq1,cmatrix2}
studying nonnegative integer solutions
to linear inequalities as they relate to the
enumeration of integer partitions and compositions.
Define the {\em weight} of a sequence
 $\lambda = (\lambda_1, \lambda_2, \ldots, \lambda_n)$
of integers
 to be
$|\lambda|=\lambda_1 + \cdots + \lambda_n$.
If sequence $\lambda$ of weight $N$
has all parts nonnegative, we call it a {\em composition} of $N$;
if, in addition, $\lambda$ is a
nonincreasing sequence, we call
it a {\em partition} of $N$.

Given an $r \times n$ integer matrix $C=[c_{i,j}]$, we consider the set $S_C$
of nonnegative integer sequences
$\lambda = (\lambda_1, \lambda_2, \ldots, \lambda_n)$ satisfying
the constraints
\begin{equation}
c_{i,0} +  c_{i,1}\lambda_1 +
c_{i,2}\lambda_2 + \ldots
+ c_{i,n}\lambda_n \geq  0, \ \ \
 1 \leq i \leq r.
\label{constraints}
\end{equation}
We seek
the full generating function
\begin{equation}
F_C(x_1,x_2, \ldots, x_n) =
\sum_{\la \in S_C}x_1^{\la_1}x_2^{\la_2} \cdots x_n^{\la_n},
\label{full}
\end{equation}
which can be viewed as an
encapsulation of the solution
set $S_C$: the coefficient of $q^N$ in $F_C(qx_1,qx_2, \ldots,
qx_n)$ is a listing (as the terms of a polynomial) of all
nonnegative integer
solutions to (\ref{constraints}) of weight $N$ and the number of
such solutions is  the coefficient of $q^N$ in $F_C(q,q, \ldots, q)$.

Variations of this problem arise in other areas of mathematics,
e.g., solving systems of linear equations,
finding volume of polytopes, as well as in enumeration.  
In the papers \cite{ineq1,cmatrix2} we demonstrated that in the area of
partition and composition enumeration many familiar sets of 
linear constraints
can be easily handled a matrix inversion:
for homogeneous systems,
if the constraint matrix $C$ is an $n \times n$ invertible
matrix, and if all entries of
$C^{-1}=B=[b_{i,j}]$ are {\em nonnegative integers}
then by Theorem 1 in \cite{cmatrix2}:
\[
F_C(x_1,x_2, \ldots, x_n) =
\prod_{j =1}^{n}
\frac{1}
 {(1-x_1^{b_{1,j}}x_2^{b_{2,j}}
 \cdots x_n^{b_{n,j}}) }.
\]
This theorem (in its full generality) suffices to handle the enumeration
of such families as
{ Hickerson partitions} \cite{hickerson},
{ Santos' interpretation of Euler's family} \cite{santos},
{ Sellers' generalization of Santos} \cite{sellers,sellers2},
{ partitions with nonnegative second differences} \cite{PA2},
{ super-concave partitions} \cite{snellman},
{ partitions with $r$-th differences nonnegative} \cite{PA2,CCH,SC},
partitions with mixed difference conditions \cite{PA2}, and
{ examples {\rm (}0-5{\rm )}  of Pak in} \cite{pak}.
The theorem provides bijections as well as generating functions.

However,
it is easy to find simple examples where the ``C matrix'' technique
fails.  In Section 2,  we propose
five simple guidelines  for computing the generating function
$F_{\C}$ of
a system $\C$ of linear diophantine
inequalities.  
The guidelines can be viewed as a simplification of
MacMahon's partition analysis \cite{Mamo2}, with ideas from
matrix
methods, Elliott reduction \cite{elliott}, and
  ``adding a slice'' (e.g. \cite{KP}).

Our focus is on the use of the guidelines to 
{\em derive a recurrence}
for the generating function $F_{\C_n}$ of an {\em infinite family}
$\{\C_n | n \geq 1\}$ of constraint systems.  
This is in contrast
to the focus of the Omega package \cite{PA3}, a software implementation of 
partition
analysis,  well-designed to compute the generating function of a given fixed,
finite system of linear constraints.
The advantage of a recurrence for  $F_{\C_n}$ is
a program which computes $F_{\C_n}$ for any given $n$.
But more significantly,
if the recurrence can be solved,  it provides a closed
form for the generating function for the infinite family.

In Sections 3-6,
we show how to use the guidelines of Section 2 strategically to solve some
nontrivial enumeration problems in the theory of partitions and
compositions.  Sections 3 and 4 address well-studied problems,
included as ``warm-up'' exercises to illustrate the approach and the
handling of the recurrences that result.  Sections 5 and 6 apply the method
to the problem of enumerating
anti-lecture hall compositions \cite{anti} and
truncated lecture hall partitions \cite{CS2},
giving a simpler approach than in \cite{anti,CS2}.
For completeness, in Section 7 we prove that the guidelines suffice to 
find the
generating function for the nonnegative integer solutions of
any homogeneous system of linear inequalities with
integer coefficients.

This work was inspired by the the work of Andrews, Paule, and
Riese in the sequence of papers 
\cite{PA1,PA2,PA3,PA4,PA5,PA6,PA7,PA8,PA9,PA10,PA11}, which illustrate many
applications of partition analysis.  The Omega Package
software \cite{PA3}
was an invaluable tool in our early investigations.
As illustrated in papers such as \cite{PA1,PA2,PA5,PA8,PA10},
recurrences can certainly be derived using partition analysis.
However, we found that the task became easier with a simpler
set of tools which appear to be no less powerful.
In Section 8 we discuss MacMahon's partition analysis
and show how the proposed guidelines can be be viewed as essential ideas
underlying his theory.

\pagebreak
\section{The Five Guidelines}
\label{section:guidelines}

Let $\C$ be a set of linear constraints in $n$ variables,
$\la_1, \ldots, \la_n$,
each constraint $c \in \C$ of the form 
\[
c:\ \  [a_0 + \sum_{i=1}^n a_i \la_i \geq 0],
\]
for integer values $a_0,a_1, \ldots, a_n$.

Let $S_{\C}$ be the set of  nonnegative
 integer sequences $ \la = (\la_1, \ldots, \la_n)$
satisfying all constraints in $\C$.
Since we are only interested here in {\em nonnegative} integer
solutions, {\em we will always assume that $\C$ contains the constraints
$[\lambda_i \geq 0]$ for $1 \leq i \leq n$}.
Define the {\em full generating function of $\C$} to be:
\[
F_{\C}(x_1, \ldots, x_n) \triangleq \sum_{\la \in S_{\C}}x_1^{\la_1}x_2^{\la_2}\cdots x_n^{\la_n}.
\]

If $c$ is the constraint:  $[a_0 + \sum_{i=1}^n a_i \la_i \geq 0]$
define the {\em negation} of $c$,
$\neg c$, to be the constraint
$[-a_0 - \sum_{i=1}^n a_i \la_i \geq 1]$.
Then any sequence $(\la_1, \ldots, \la_n)$ satisfies
$c$ or $\neg c$, but not both.
A constraint $c$ is {\em implied by} the set of constraints $\C$ if
$S_{\C \cup \{\neg c\}}= \emptyset$.
A constraint $c$ is {\em redundant} if
$S_{\C \cup \{ c\}}= S_{\C}$.

Let
$\C_{\la_i \leftarrow \la_i+a\la_j}$ denote the set of constraints
which results  from  replacing  $\la_i$ by
 $\la_i+a\la_j$ in every constraint in $\C$.
Note that if constraint $c$ is implied by $\C$ then
$c_{\la_i \leftarrow \la_i+a\la_j}$ is implied by
$\C_{\la_i \leftarrow \la_i+a\la_j}$.
Thus observe that if $\C$ contains the constraints
$[\la_k \geq 0], 1 \leq k \leq n$ and if
$[\la_i - a\la_j \geq 0]$ is implied by $\C$, then
all of the constraints $[\la_k \geq 0], 1 \leq k \leq n$
are also implied by 
$\C_{\la_i \leftarrow \la_i+a\la_j}$.

\begin{lemma}
Let $\C$ be a set of linear constraints on variables
$\la_1, \ldots, \la_n$ which
contains the constraints
$[\la_k \geq 0], 1 \leq k \leq n$.
Let $a$ be any integer (possibly negative).
Suppose
$[\la_i-a\la_j\geq 0]$
is implied by  $\C$ and let
$\C' = \C_{\la_i \leftarrow \la_i+a\la_j}$. Then
\[
\beta = (\beta_1, \ldots, \beta_n) \in S_{\C} \ \ \   {\rm iff} \ \ \ 
\beta' = (\beta_1, \ldots, \beta_{i-1}, \beta_i-a\beta_{j},
 \beta_{i+1}, \ldots,  \beta_n)
 \in S_{\C'}.
\]
\label{columnops}
\end{lemma}
\noindent
{\bf Proof.}
By the remarks preceding the lemma, the constraints $\C$ and $\C'$
guarantee that $ S_{\C}$ and  $ S_{\C'}$
contain only nonnegative integer solutions.  So, it suffices to show
that  $\beta$ satisfies a constraint in $\C$ iff
$\beta'$ satisfies the corresponding constraint in $\C'$.

Let $c(\la) = c_0 + \sum_{t=1}^n c_t \la_t$ and assume
$[c(\la) \geq 0] \in \C$.
Under the substitution
$\la_i \leftarrow \la_i+a\la_j$, $c(\la)$ becomes $c'(\la)$
defined by
\[
c'(\la) = c_0 + \sum_{t=1}^n c_t \la_t + c_ia\la_j = c(\la) + c_ia\la_j
\]
and $[c'(\la) \geq 0] \in \C'$.  Thus
$$c(\beta) = c'(\beta) - c_ia\beta_j =c'(\beta'),$$
so $c(\beta) \geq 0$ iff $c'(\beta') \geq 0$.
\qed

%We first check that $\beta$ satisfies a constraint in $\C$ iff
%$\beta'$ satisfies the corresponding constraint in $\C'$.
%Let $c(\la) = c_0 + \sum_{t=1}^n c_t \la_t$ and assume
%$[c(\la) \geq 0] \in \C$.
%Under the substitution
%$\la_i \leftarrow \la_i+a\la_j$, $c(\la)$ becomes $c'(\la)$
%defined by
%\[
%c'(\la) = c_0 + \sum_{t=1}^n c_t \la_t + c_ia\la_j = c(\la) + c_ia\la_j
%\]
%and $[c'(\la) \geq 0] \in \C'$.  Thus
%\[
%c(\beta) = c'(\beta) - c_ia\beta_j =c'(\beta'),
%\]
%so $c(\beta) \geq 0$ iff $c'(\beta') \geq 0$.

%Now check that if $\beta \in S_{\C}$ then $\beta'$ is a nonnegative
%integer sequence.
%By hypothesis,
%$[\la_i-a\la_j\geq 0]$
%is implied by $\C$ and $\C$  
%contains the constraints
%$[\la_i \geq 0], 1 \leq i \leq n$.
%So $\beta \in S_{\C}$ implies 
%$\beta_i-a\beta_j\geq 0$ and
%$\beta_i\geq 0, 1 \leq i \leq n$. 
%Conversely, since 
%$\C$ contains $[\la_i \geq 0], 1 \leq i \leq n$,
%$\C'$ contains the constraint $[\la_i + a \la_j \geq 0]$ and
%the constraints
% $[\la_t \geq 0], 1 \leq t < i, i < t  \leq n$.
%Thus,
% $\beta' \in S_{\C'}$ implies
%$\beta'_i + a \beta'_j \geq 0$, i.e,
%$(\beta_i-a\beta_j)+a\beta_j \geq 0$ and
% $[\beta_t \geq 0], 1 \leq t < i, i < t  \leq n$,
%so $\beta$ is nonnegative.
%\qed

Finally, to simplify notation, we will let $X_n$ refer to the parameter
list $x_1, \ldots, x_n$, so that $F(X_n)$ denotes
$F(x_1, \ldots, x_n)$.
Let $F(X_n;x_i \leftarrow x_ix_j^a)$ denote the
function $F(X_n)$ with all  occurrences of $x_i$  replaced
by $x_ix_j^a$.

\begin{theorem}
(The Five Guidelines) 

 1. If $\C=\{[\la_1 \geq t]\}$, for integer $t \geq 0$, then
\[
F_{\C}(x_1) = \frac{x_1^t}{1-x_1}.
\]

2. If $\C_1$ is a set of constraints on variables
$\la_1, \ldots, \la_j$ and $\C_2$
is a set of constraints on variables $\la_{j+1}, \ldots, \la_n$, then
\[
F_{\C_1 \cup \C_2}(x_1, \ldots, x_n)=F_{\C_1}(x_1, \ldots x_j)
F_{\C_2}(x_{j+1}, \ldots, x_n).
\]

3.   
Let $\C$ be a set of linear constraints on variables
$\la_1, \ldots, \la_n$ and assume
$\C$ contains the constraints
$[\la_i \geq 0], 1 \leq i \leq n$.
Let $a$ be any integer (possibly negative).
If $[\la_i-a\la_j\geq 0]$
 is implied by  $\C$,
\[
F_{\C}(X_n)=F_{\C_{\la_i \leftarrow \la_i+a\la_j}}(X_n;x_j \leftarrow x_jx_i^a).
\]

4. Let $c$ be any constraint with the same variables as the set
$\C$.  Then
\[
F_{\C}(X_n)=F_{\C \cup \{c\}}(X_n) + F_{\C \cup \{\neg c\}}(X_n).
\]

5. Let $c \in \C$.  Then
\[
F_{\C}(X_n)=F_{\C-\{c\}}(X_n)-F_{\C-\{c\}\cup\{\neg c\}}(X_n).
\]
\label{thm:guidelines}
\end{theorem}
\noindent
{\bf Proof.}

1. This is clear since $F_{\C}(x_1) =x_1^t + x_1^{t+1} + \cdots$.

2. The sequence $(\la_1, \ldots, \la_n) \in S_{\C_1 \cup \C_2}$ iff
$(\la_1, \ldots, \la_j) \in S_{\C_1 }$  and
 $(\la_{j+1}, \ldots, \la_n) \in S_{\C_2}$.

3. Let $\C' = \C_{\la_i \leftarrow \la_i+a\la_j}$.
By Lemma \ref{columnops}, 
\[
(\la_1, \ldots, \la_n) \in S_{\C'} \ \ \   {\rm iff} \ \ \ 
(\la_1, \ldots, \la_{i-1}, \la_i+a\la_{j}, \la_{i+1}, \ldots,  \la_n)
 \in S_{\C}.
\]
So,
\begin{eqnarray*}
F_{\C'}(X_n;x_j \leftarrow x_jx_i^a) & = &
\sum_{\la \in S_{\C'}}x_1^{\la_1}x_2^{\la_2}\cdots
x_{j-1}^{\la_{j-1}}
(x_{j}x_i^a)^{\la_{j}}
x_{j+1}^{\la_{j+1}}
\cdots
 x_n^{\la_n} \\
& = &
\sum_{\la \in S_{\C'}}x_1^{\la_1}x_2^{\la_2}\cdots
x_{i-1}^{\la_{i-1}}
x_{i}^{(\la_{i}+a\la_j)}
x_{i+1}^{\la_{i+1}}
\cdots
 x_n^{\la_n} \\
& = &
\sum_{\la \in S_{\C}}x_1^{\la_1}x_2^{\la_2}\cdots
x_{i}^{\la_{i}}
\cdots
 x_n^{\la_n} \\
& = & F_{\C}(X_n).
\end{eqnarray*}

4.  $S_{\C}$ can be partitioned into those $\la$
that satisfy $c$ and those that do not.

5. By guideline 4,
$F_{\C-\{c\}}(X_n)= 
F_{\C - \{c\} \cup \{c\}}(X_n) + F_{\C - \{c\} \cup \{\neg c\}}(X_n)$.
%Since $c \in \C$, $\C - \{c\} \cup \{c\}=C$.
Then 
$\C - \{c\} \cup \{c\}=\C$,
since $c \in \C$. 
\qed

\pagebreak
\section{Minc's Partition Function and Cayley Compositions}

Minc's partition function  $\nu(d,N)$ is
the number of compositions of $N$ in which the first part is $d$ and each part
is at most twice the size of the preceding part \cite{minc}.
For example, 
in the special case $d=1$, these are called Cayley compositions
\cite{cayley,RRreciprocal,PA5}.  In this section we compute the generating
function $ \nu(q) = \sum_{d,N \geq 0} \nu(d,N) q^N=
q+2q^2+4q^3+7q^4+13q^5+24q^6+ \cdots$.
For example, the coefficient of $q^5$ is 13, since of the 16 compositions
of 5, only these three violate the constraints:
$(1,4)$, $(1,3,1)$, and $(1,1,3)$.

Let ${\mathcal C}_n$ be the set of constraints
${\mathcal C}_n  =  \{\la_i \geq \frac{1}{2} \la_{i+1} >0 \ | \ 1 \leq i < n\}$
and let $C_n(x_1, \ldots, x_n)$ be the generating function of 
${\mathcal C}_n$.
Focusing on the constraint
$c=[\la_{n-1}  \geq  \frac{1}{2} \la_{n}]$,
after noting that $[\la_{n-1}>0]$ is redundant,
we can write ${\mathcal C}_n$ as
{\small
\[
{\mathcal C}_n =
\left [ \begin{array}{rcl}
\la_1 & \geq  &\frac{1}{2} \la_2 \\
\vspace{-.05in}\\
\la_2  &\geq  &\frac{1}{2} \la_3 \\
&\vdots &\\
\la_{n-2}  &\geq  & \frac{1}{2} \la_{n-1} \\
\vspace{-.05in}\\
\la_{n-1}  &\geq  & \frac{1}{2} \la_{n} \\
\vspace{-.05in}\\
\la_n &  > &0
\end{array} \right ]
 \  = \
\left [ \begin{array}{rcl}
\la_1 & \geq  &\frac{1}{2} \la_2 \\
\vspace{-.05in}\\
\la_2  &\geq  &\frac{1}{2} \la_3 \\
&\vdots &\\
\la_{n-2}  &\geq  & \frac{1}{2} \la_{n-1} \\
\vspace{-.05in}\\
\la_{n-1}  &\geq  & \frac{1}{2} \la_{n} \\
\vspace{-.05in}\\
\la_{n-1}  & > & 0  \\
\vspace{-.05in}\\
\la_n &  > &0
\end{array} \right ]
 \  = \
\left [ \begin{array}{rcl}
\la_1 & \geq  &\frac{1}{2} \la_2 \\
\vspace{-.05in}\\
\la_2  &\geq  &\frac{1}{2} \la_3 \\
&\vdots &\\
\la_{n-2}  &\geq  & \frac{1}{2} \la_{n-1} \\
\vspace{-.05in}\\
\la_{n-1}  &>  & 0 \\
\vspace{-.05in}\\
\la_{n}  & > &0
\end{array} \right ]
 \  - \
\left [ \begin{array}{rcl}
\la_1 & \geq  &\frac{1}{2} \la_2 \\
\vspace{-.05in}\\
\la_2  &\geq  &\frac{1}{2} \la_3 \\
&\vdots &\\
\la_{n-2}  &\geq  & \frac{1}{2} \la_{n-1} \\
\vspace{-.05in}\\
\la_{n}  &>  & 2\la_{n-1} \\
\vspace{-.05in}\\
\la_{n-1}  &>  & 0 
\end{array} \right ],
\]
}
where $c$ has been removed from the next-to-last system,
making it ${\mathcal C}_{n-1} \cup [\la_n>0]$, and $c$ has been
replaced by $\neg c$ in the last system.
By guidelines 1 and 2,
$x_n C_{n-1}(x_1, \ldots, x_{n-1})/(1-x_n)$ is the generating function
for
${\mathcal C}_{n-1} \cup [\la_n>0]$.
Note further that the
substitution $\la_n \leftarrow \la_n+2\la_{n-1}$ in the last system
results ${\mathcal C}_{n-1} \cup [\la_n>0]$, so by guideline 3,
the last system has generating function
$x_n C_{n-1}(x_1, \ldots, x_{n-1}x_n^2)/(1-x_n)$. 
Putting this together with guideline 5 and the initial condition
$C_1(x_1)=x_1/(1-x_1)$ 
gives the
recurrence
\[
C_n(x_1,\ldots, x_n) = \frac{x_n}{1-x_n}(C_{n-1}(x_1, \ldots, x_{n-1})-
C_{n-1}(x_1, \ldots, x_{n-2},x_{n-1}x_n^2)).
\]
Let $C_n(q,s)=C_n(q,q,\ldots,q,s)$. Then the above
recurrence gives $C_1(q,s)=s/(1-s)$ and for $n \geq 2$,
\[
C_n(q,s) = \frac{s}{1-s}(C_{n-1}(q,q)-C_{n-1}(q,qs^2)).
\]
Set $C(q,s)\ =\ \sum_{n=1}^{\infty}C_n(q,s)$ and use the recurrence for 
$C_n(q,s)$ to get
\[
C(q,s)\ =\ \sum_{n=1}^{\infty}C_n(q,s)\ =\ 
\frac{s}{1-s}+\sum_{n=2}^{\infty}C_n(q,s)
\ = \  \frac{s}{1-s}(1+C(q,q) - C(q,qs^2)).
\]
Iterating the recurrence for $C(q,s)$ gives
\[
C(q,s)=(1+C(q,q))\sum_{i=1}^{\infty}
(-1)^{i-1}
\prod_{j=0}^{i-1}\frac{q^{2^j-1}s^{2^j}}{(1-q^{2^j-1}s^{2^j})}.
\]
Let $C(q)=C(q,q)$, then
\[
\nu(q)  =
1+C(q)= \frac{1}{1+ \sum_{i=1}^{\infty}
\frac{(-1)^i q^{2^{i+1}-i-2}}{(1-q)(1-q^3)(1-q^7) \cdots (1-q^{2^{i}-1})}}.
\]
%(This ends up being the same as Andrews without the $q$ in the
%numerator, but his generating function is wrong.  I could go ahead
%and do the calculation for first part = c (or d).)
%
%Let $c_j(n)$ be the number of Cayley compositions with $j$ positive parts
%and last part $n$.  Let
%\[
%C_j(q)= \sum_{n=0}^{\infty}c_j(n)q^n.
%\]
%Then $C_j(q)=p_{j-1}(1,1, \ldots,1,q)$.  So the recurrence for $p$
%gives
%\[
%C_j(q)=\frac{q}{1-q}C_{j-1}(1)-C_{j-1}(q^2).
%\]
%
%{\em Cayley's Theorem} says that the number of Cayley compositions
%with $j$ parts is the number of partitions of $2^{j-1}-1$ into
%parts in $\{1,1',2,4,8, \ldots, 2^{j-2}\}$.
%
%(So note this is saying something about $C_j(1)$.)
%
%To prove this, they get the alternating sum recurrence for $C_j(q)$
%and then do something clever by looking also at
%$C_j(1/q)$.

%Seems like this could be easier?
%
%Anyway - notes:
%
%- The infinite form of the generating function for $p_n$ is
%found in \cite{RRreciprocal} - see Section \ref{section:reciprocal}.
%(It turns out Minc wanted to count these in \cite{minc}).
%
%- Everybody knows it is easy to generalize this to $1/k$ instead of
%$1/2$.
%
%- In the $1/2$ case, what if we want $\la_1=c$?  (Minc considered these
%even though he really only wanted  $\la_1=1$ because he used it
%to get his recurrence.)

%- We might think of anti-lecture hall compositions as a generalization
%of this and ask how many with first part 1? See Section \ref{section:cayley-anti}.
%

\pagebreak
\section{Two-Rowed Plane Partitions}

This example illustrates the advantage of 
guideline 3 of Theorem \ref{thm:guidelines} when
$a<0$.
The {\em two-rowed plane partitions} are those  integer
sequences $(a_1,b_1, \ldots, a_n,b_n)$ satisfying the constraints
\[
\mP_n  \ \ \ = \ \  \   
\left [ a_i \geq b_i \geq 0,  \ \  1 \leq i \leq n;
 \ \ \ \ \ \ \ \  a_i \geq a_{i+1}, \ \ \  b_i \geq b_{i+1},
  \ \ \  1 \leq i \leq n-1 \right ].
\]
It is well-known that the generating function for $\mP_n$ is
\cite{MacMahon-ppgf}
\begin{equation}
P_n(q) = \frac{1}{(q;q)_n(q^2;q)_n}.
\label{2row}
\end{equation}
In \cite{PA2}, Andrews shows how MacMahon's {\em partition analysis}
can be used to
compute  $P_n(q)$ by considering an
intermediate family $\G_n$.
We will use this approach, but with a slight twist,
to show how the generating function for
 $\mP_n$, can be computed via  $\G_n$ from
the guidelines of Theorem \ref{thm:guidelines}.

%For a constraint system $\C$, by {\em generating function of $\C$} we
%mean the generating function of $S_{\C}$.
We will use the convention that when a constraint system is
represented by a calligraphic letter, its generating function
is represented by  the corresponding roman letter.
Also, to keep notation simple, when the meaning is clear from
context, we will use the same letter to
refer to multivariable and single variable forms of the
generating function.

Define  $\G_n$ to be the set of constraints below:
{\small
\[
\G_n  \ \ \ = \ \ \
\left [ \begin{array}{rcr}
a_1 + a_2 + \cdots +a_n & \geq & b_1 + b_2 + \cdots +b_n\\
 a_2 + \cdots +a_n & \geq &  b_2 + \cdots +b_n\\
 \vdots   & \vdots  & \vdots \\
a_{n-1}+a_n & \geq & b_{n-1}+b_n\\
a_{n} & \geq & b_n\\
a_i, b_i \geq 0, & &  i=1, \ldots, n
\end{array} \right ].
\]
}

\noindent
Denote the full generating functions for $\mP_n$ and $\G_n$ by
\[
P_n(x_1,y_1, \ldots, x_n,y_n) \triangleq
\sum_{(a_1,b_1, \ldots, a_n,b_n) \in S_{\mP_n}}
x_1^{a_1}y_1^{b_1} \ldots, x_n ^{a_n}y_n^{b_n}, 
\]
\[
G_n(x_1,y_1, \ldots, x_n,y_n) \triangleq
\sum_{(a_1,b_1, \ldots, a_n,b_n) \in S_{\G_n}}
x_1^{a_1}y_1^{b_1} \ldots, x_n ^{a_n}y_n^{b_n}. 
\]

\noindent
Note that  $\mP_n$ can be transformed into $\G_n$
by the sequence of substitutions:
\[
a_i \leftarrow a_i+a_{i+1}; \ \ \ \ \ \ \  b_i \leftarrow b_i
+b_{i+1}; \ \ \ \ \ \ \ \   i = 1,2, \ldots n-1.
\]

\noindent
We focus on $G_n$.
Since  for $1 \leq i \leq n-1$,
$a_i -a_{i+1} \geq 0$ and
$b_i -b_{i+1} \geq 0$ in $\mP$,
by guideline 3 of Theorem \ref{thm:guidelines},
$P_n$ is obtained from  $G_n$ by the sequence of
substitutions:
\[
x_i \leftarrow x_ix_{i-1};\ \ \ \ \ \ \ y_i \leftarrow y_i y_{i-1}
\ \ \ \ \ \ \ \ i = n, n-1, n-2, \ldots, 2.
\]

\noindent
Thus
\[
P_n(x_1,y_1, \ldots, x_n,y_n) = G_n(x_1,y_1, x_1x_2,y_1y_2, \ldots,x_1x_2
\cdots x_n, y_1y_2 \cdots y_n).
\]

\noindent
In particular, the generating function (\ref{2row}) for two-rowed plane
partitions is obtained by setting $x_i=y_i=q$ in $P_n$ 
for $i = 1, \ldots, n$: 
\begin{equation}
P_n(q,q,q,\ldots,q)=G_n(q,q,q^2,q^2, \ldots, q^n,q^n).
\label{Pqqq}
\end{equation}

\noindent
Since $a_n -b_n \geq 0$ in $\G_n$, by guideline 3, 
we can do the substitution
$a_n \leftarrow a_n+b_n$ in $\G_n$ to get $\F_n$ 
and recover  $G_n$ from $F_n$ as
shown below.  
{\small
\[
\F_n \ \ \ \ 
 = \ \ \ 
\left [
\begin{array}{rcr}
a_1 + a_2 + \cdots +a_n & \geq & b_1 + b_2 + \cdots +b_{n-1}\\
 a_2 + \cdots +a_n & \geq &  b_2 + \cdots +b_{n-1}\\
  \vdots  & \vdots  & \vdots \\
a_{n-1}+a_n & \geq & b_{n-1}\\
a_i, b_i \geq 0, & &  i=1, \ldots, n
\end{array}
\right ],
\]}
\[
G_n(x_1,y_1, \ldots, x_n,y_n) =
F_n(x_1,y_1, \ldots, x_n,y_n; y_n \leftarrow x_ny_n ).
\]

\noindent
Since $a_{n-1}+a_n \geq 0$ in $\F_n$, by guideline 3, 
we can substitute
$a_{n-1} \leftarrow a_{n-1}-a_n$ in $\F_n$ to get $\mH_n$ 
and recover $F_n$ from $H_n$ as shown.
{\small
\[
\mH_n  \ \ \ \  = \ \ \ \
\left [
\begin{array}{rcr}
a_1 + a_2 + \cdots +a_{n-1} & \geq & b_1 + b_2 + \cdots +b_{n-1}\\
 a_2 + \cdots +a_{n-1} & \geq &  b_2 + \cdots +b_{n-1}\\
   \vdots & \vdots  & \vdots \\
a_{n-1} & \geq & b_{n-1}\\
a_{n-1} & \geq & a_n\\
a_i, b_i \geq 0 & & i = 1, \ldots, n
\end{array} \right ],
\]
}
\[
F_n(x_1,y_1, \ldots, x_n,y_n) =
H_n(x_1,y_1, \ldots, x_n,y_n; x_n \leftarrow x_n/x_{n-1} ).
\]
Summarizing to this point, we have
\begin{equation}
G_n(x_1,y_1, \ldots, x_n,y_n) =
H_n(x_1,y_1, \ldots,x_{n-1},y_{n-1}, x_n/x_{n-1},x_ny_n).
\label{Gsofar}
\end{equation}

\noindent
Now apply guideline 5 to $\mH_n$ using the constraint
$c = [a_{n-1} \geq a_n]$.  Then  
\[
\mH_n=\K_n-\mL_n, 
\]
where
\[
\K_n =\mH_n - \{[a_{n-1} \geq a_n]\},  \ \ \ \ \ 
\mL_n = \mH_n - \{[a_{n-1} \geq a_n]\} \cup \{[a_{n} \geq  a_{n-1}+1]\},
\]
that is,
\[
\K_n \ \ \  =  \ \ \ 
\left [
\begin{array}{rcr}
a_1 + a_2 + \cdots +a_{n-1} & \geq & b_1 + b_2 + \cdots +b_{n-1}\\
 a_2 + \cdots +a_{n-1} & \geq &  b_2 + \cdots +b_{n-1}\\
   \vdots & \vdots  & \vdots \\
a_{n-1} & \geq & b_{n-1}\\
a_i, b_i \geq 0 & & i = 1, \ldots, n
\end{array} \right ]
\]
and
\[
\mL_n \ \ \  =  \ \ \ 
\left [
\begin{array}{rcr}
a_1 + a_2 + \cdots +a_{n-1} & \geq & b_1 + b_2 + \cdots +b_{n-1}\\
 a_2 + \cdots +a_{n-1} & \geq &  b_2 + \cdots +b_{n-1}\\
   \vdots & \vdots  & \vdots \\
a_{n-1} & \geq & b_{n-1}\\
a_n & \geq  & a_{n-1}+1\\
a_i, b_i \geq 0 & & i = 1, \ldots, n
\end{array} \right ],
\]
 so that
\begin{equation}
H_n(x_1,y_1, \ldots, x_n,y_n) =
K_n(x_1,y_1, \ldots, x_n,y_n) -
L_n(x_1,y_1, \ldots, x_n,y_n).
\label{KminusL}
\end{equation}

\noindent
Now observe that
\[
\K_n = \G_{n-1} \cup \{[a_n \geq 0], [b_n \geq 0]\},
\]
so by guidelines 1 and 2,
\begin{equation}
K_n(x_1,y_1, \ldots, x_n,y_n)=
\frac{G_{n-1}(x_1,y_1, \ldots, x_{n-1},y_{n-1})}{(1-x_n)(1-y_n)}.
\label{Kn}
\end{equation}

\noindent
Returning to $\mL_n$,
since $a_n-a_{n-1} \geq 0$ in $\mL_n$, we
can do the substitution
$a_n \leftarrow a_n+a_{n-1}$, resulting in
\[
(\mL_n)_{a_n \leftarrow a_n+a_{n-1}} = 
\G_{n-1} \cup \{[a_n \geq 1], [b_n \geq 0]\},
\]

\noindent
so by guidelines 1, 2, and 3,
\begin{equation}
L_n(x_1,y_1, \ldots, x_n,y_n)=
\frac{x_n G_{n-1}(x_1,y_1, \ldots, x_{n-1},y_{n-1}; x_{n-1} \leftarrow
x_{n-1}x_n)}{(1-x_n)(1-y_n)}.
\label{Ln}
\end{equation}

\noindent
Combining (\ref{KminusL}),(\ref{Kn}), and (\ref{Ln}),
we have
\begin{equation*}
H_n(x_1,y_1, \ldots, x_n,y_n)=\frac{G_{n-1}(x_1,y_1, \ldots, x_{n-1},y_{n-1})}
{(1-x_n)(1-y_n)}-\frac{x_nG_{n-1}(x_1,y_1, \ldots, x_{n-2},y_{n-2},x_{n-1}x_n,y_
{n-1})}
{(1-x_n)(1-y_n)}.
\end{equation*}

\noindent
Finally,
substituting this expression for $H_n$ into (\ref{Gsofar})
gives a recurrence for $G_n$:
\begin{equation}
G_n(x_1,y_1, \ldots, x_n,y_n)=\frac{G_{n-1}(x_1,y_1, \ldots, x_{n-1},y_{n-1})
-\frac{x_{n}}{x_{n-1}}G_{n-1}(x_1,y_1, \ldots, x_{n
-2},y_{n-2},x_n,y_{n-
1})}
{(1-x_{n}/x_{n-1})(1-x_ny_n)},
\label{Grec}
\end{equation}

\noindent
with initial condition $G_1(x_1,y_1)=1/(1-x_1)/(1-x_1y_1)$.

Let $G_n^*(q,s)=G_n(q,q,q^2,q^2, \ldots,s,q^n)$.  Then from the
recursion (\ref{Grec}),
\[
G_n^*(q,s)=\frac{G^*_{n-1}(q,q^{n-1})-(s/q^{n-1})
G^*_{n-1}(q,s)}{(1-s/q^{n-1})(1-sq^n)}.
\]

\noindent
It is straightforward to show by induction that $G_n^*(q,s)$ satisfies
\[
G_n^*(q,s)= \frac{1}{(1-s)(1-sq)(q;q)_{n-1}(q^2;q)_{n-1}}.
\]

\noindent
Substituting $s=q^n$ gives
\[
P_n(q) = 
G_n(q,q,q^2,q^2, \ldots,q^n,q^n)=G_n^*(q,q^n)=\frac{1}{(q;q)_n(q^2;q)_n},
\]
the desired generating function for $2 \times n$ plane partitions.

\pagebreak
\section{Anti-Lecture Hall Compositions}
\label{section:anti}

In \cite{anti}, we considered
the
 set of sequences $\la=(\la_1, \ldots, \la_n)$
satisfying the constraints
\[
\A_n \ \ \ = \ \ \ \left [ 
\frac{\la_1}{1}\ge\frac{\la_2}{2}\ge \ldots \ge \frac{\la_n}{n}\ge 0
\right ].
\]
We referred to these as
{\em anti-lecture hall compositions}
and showed that the generating function
is
\begin{equation}
A_n(q) \triangleq
\sum_{\lambda\in A_n} q^{|\lambda|} =
\prod_{i=1}^{n}\frac{1+q^i}{1-q^{i+1}}.
\label{antiq}
\end{equation}

Here we show how to apply the guidelines of Theorem \ref{thm:guidelines}
to get a recurrence for the full generating function
$A_{n}(x_1,x_2, \ldots x_n)$ and use it to give
an ``easy'' proof of
(\ref{antiq}).
The idea is easily extended
to the  {\em truncated anti-lecture hall compositions}
studied in \cite{CS2}.
We start with $\B_n$, a slight variation of $\A_n$.
\begin{lemma}
The full generating function for the integer
sequences defined by the constraints
\begin{equation}
\B_n \ \ \ = \ \ \ \left [
\frac{\la_1}{1} \geq
\frac{\la_2}{2} \geq
\cdots \geq
\frac{\la_{n-1}}{n-1} \geq
\frac{\la_n}{1} \geq 0 \right].
\label{Bconstraints}
\end{equation}
satisfies
\[
B_n(x_1,\ldots,x_n)=
\frac{A_{n-1}(x_1 \ldots,x_{n-1})}{1-x_1x_2^2x_3^3 \cdots x_{n-1}^{n-1}x_n}.
\]
\label{lemma:Bn}
\end{lemma}
\noindent
{\bf Proof.}
The following sequence of substitutions transforms $\B_n$ into
$\A_{n-1} \cup \{[\la_n \geq 0]\}$, as illustrated in 
Figure \ref{figure:anti}:
\[
\la_i \leftarrow \la_i+i\la_n, \ \ \  i= n-1, \ldots, 1.
\]
Note that the constraint
$\la_{i-1} \geq (i-1)\la_n$ is implied at each stage, so by guidelines
1,2, and 3, $B_n$ is recovered from $A_n$ by performing the
sequence of substitutions on $A_{n-1}(x_1, \ldots, x_{n-1})/(1-x_n)$:
\[
x_n \leftarrow  x_n x_{i}^{i}, \ \ \ 
i =  1, \ldots, n-1.
\]
\qed
\hspace{-.5in}
\begin{figure}
{\small
\[
\left [ \begin{array}{rcl}
\la_1 & \geq  &\frac{1}{2} \la_2 \\ 
\vspace{-.05in}\\
\la_2  &\geq  &\frac{2}{3} \la_3 \\
&\vdots &\\
\la_{n-3} & \geq  &\frac{n-3}{n-2} \la_{n-2} \\
\vspace{-.05in}\\
\la_{n-2} & \geq  &\frac{n-2}{n-1} \la_{n-1} \\
\vspace{-.05in}\\
\la_{n-1}  &\geq  &(n-1) \la_{n} \\
\vspace{-.05in}\\
\la_n  &\geq  &0
\end{array} \right ]
 \ \rightarrow  \ 
\left [ \begin{array}{rcl}
\la_1 & \geq  &\frac{1}{2} \la_2 \\
\vspace{-.05in}\\
\la_2  &\geq  &\frac{2}{3} \la_3 \\
&\vdots &\\
\la_{n-3} & \geq  &\frac{n-3}{n-2} \la_{n-2} \\
\vspace{-.05in}\\
\la_{n-2} & \geq  &\frac{n-2}{n-1} \la_{n-1} +(n-2)\la_{n}\\
\vspace{-.05in}\\
\la_{n-1}  &\geq & 0\\
\vspace{-.05in}\\
\la_n  &\geq  &0
\end{array} \right ]
 \ \rightarrow \ 
\left [ \begin{array}{rcl}
\la_1 & \geq  &\frac{1}{2} \la_2 \\
\vspace{-.05in}\\
\la_2  &\geq  &\frac{2}{3} \la_3 \\
&\vdots &\\
\la_{n-3} & \geq  &\frac{n-3}{n-2} \la_{n-2} +(n-3)\la_{n} \\
\vspace{-.05in}\\
\la_{n-2} & \geq  &\frac{n-2}{n-1} \la_{n-1} \\
\vspace{-.05in}\\
\la_{n-1}  &\geq & 0\\
\vspace{-.05in}\\
\la_n  &\geq  &0
\end{array} \right ]
 \ \rightarrow \  
\]

\[
\rightarrow \ \ \ \ldots \ \ \ \rightarrow 
\left [ \begin{array}{rcl}
\la_1 & \geq  &\frac{1}{2} \la_2 \\
\vspace{-.05in}\\
\la_2  &\geq  &\frac{2}{3} \la_3 + 2\la_{n} \\
&\vdots &\\
\la_{n-3} & \geq  &\frac{n-3}{n-2} \la_{n-2} \\
\vspace{-.05in}\\
\la_{n-2} & \geq  &\frac{n-2}{n-1} \la_{n-1} \\
\vspace{-.05in}\\
\la_{n-1}  &\geq  & 0 \\
\vspace{-.05in}\\
\la_n  &\geq  &0
\end{array} \right ]
\ \ \ \rightarrow \ \ \
\left [ \begin{array}{rcl}
\la_1 & \geq  &\frac{1}{2} \la_2 + \la_{n} \\
\vspace{-.05in}\\
\la_2  &\geq  &\frac{2}{3} \la_3  \\
&\vdots &\\
\la_{n-3} & \geq  &\frac{n-3}{n-2} \la_{n-2} \\
\vspace{-.05in}\\
\la_{n-2} & \geq  &\frac{n-2}{n-1} \la_{n-1} \\
\vspace{-.05in}\\
\la_{n-1}  &\geq  &0 \\
\vspace{-.05in}\\
\la_n  &\geq  &0
\end{array} \right ]
\ \ \ \rightarrow \ \ \
\left [ \begin{array}{rcl}
\la_1 & \geq  &\frac{1}{2} \la_2 \\
\vspace{-.05in}\\
\la_2  &\geq  &\frac{2}{3} \la_3 \\
&\vdots &\\
\la_{n-3} & \geq  &\frac{n-3}{n-2} \la_{n-2} \\
\vspace{-.05in}\\
\la_{n-2} & \geq  &\frac{n-2}{n-1} \la_{n-1} \\
\vspace{-.05in}\\
\la_{n-1}  &\geq  &0 \\
\vspace{-.05in}\\
\la_n  &\geq  &0
\end{array} \right ]
\]
}
\caption{Transformation of $\B_n$ into $\A_{n-1} \cup \{[\la_n \geq 0]\}$
 in proof of 
Lemma \ref{lemma:Bn}.}
\label{figure:anti}
\end{figure}

\begin{proposition}
The full generating function for anti-lecture hall compositions satisfies:
\begin{eqnarray*}
A_{n}(x_1, \ldots x_n)   & = &
        \frac{A_{n-1}(x_1, \ldots,x_{n-1})}{1-x_n}
         -  {A_{n-1}(x_1, \ldots, x_{n-2},x_nx_{n-1})}\left(
         \frac{1}{1-x_n}
         -
        \frac{1}{1- x_1x_2^2x_3^3 \cdots  x_n^n} \right)
\end{eqnarray*}
with initial condition
$A_{1}(x_1) = 1/(1-x_1)$.
\label{thm:alhc-rec}
\end{proposition}
\noindent{\bf Proof.}
Using
guideline 5 with $c=[\la_{n-1} \geq \frac{n-1}{n}\la_n]$,
\[
A_n(x_1, \ldots, x_n)=
C_n(x_1, \ldots, x_n)-
D_n(x_1, \ldots, x_n),
\]
where
\begin{equation}
\C_n \ \ \ = \ \ \ \left [
\frac{\la_1}{1}\ge\frac{\la_2}{2}\ge \ldots \ge \frac{\la_{n-1}}{n-1}\ge 0;
\ \ \ \ \ \la_n \geq 0
\right ];
\label{Cn}
\end{equation}
\begin{equation}
\D_n \ \ \ = \ \ \ \left [
\frac{\la_1}{1}\ge\frac{\la_2}{2}\ge \ldots \ge \frac{\la_{n-1}}{n-1}\ge 0;
\ \ \ \ \ \frac{\la_n}{n} >  \frac{\la_{n-1}}{n-1}
\right ].
\label{Dn}
\end{equation}
Note that $\C_n = \A_{n-1} \cup \{[\la_n \geq 0] \}$, so by guideline 2,
$\C_n$ has generating function
\begin{equation}
C_n(x_1 \ldots,x_{n}) =
\frac{A_{n-1}(x_1 \ldots,x_{n-1})}{1-x_n}.
\label{first-part-gf}
\end{equation}
Since $\la_n \geq \la_{n-1}$ is implied by $\D_n$ in (\ref{Dn}),
by guideline 3,
substituting $\la_n \leftarrow \la_n+\la_{n-1}$ in $\D_n$ gives
\begin{equation}
\E_n \ \ \ = \ \ \ \left [
\frac{\la_1}{1}\ge\frac{\la_2}{2}\ge \ldots \ge \frac{\la_{n-1}}{n-1}\ge 0;
\ \ \ \ \ \la_n >  \frac{\la_{n-1}}{n-1}
\right ]
\label{En}
\end{equation}
and
\[
D_n(x_1, \ldots, x_n) = E_n(X_n; x_{n-1} \leftarrow x_{n-1} x_{n}),  
\]
where $X_n$ represents the argument list $x_1, \ldots, x_n$.
Using
guideline 5 again, with $c=[\la_n > \frac{\la_{n-1}}{n-1}]$,
gives
\[
E_n(X_n)=
C_n(X_n)-
B_n(X_n),
\]
where $\C_n$ is (\ref{Cn}) and
where $\B_n$ is (\ref{Bconstraints}).
Putting this all together, we have
\begin{eqnarray*}
A_n(X_n) & = & C_n(X_n)- D_n(X_n)\\
         & = & C_n(X_n)- E_n(X_n; x_{n-1} \leftarrow x_{n-1} x_{n})\\
         & = & C_n(X_n)- C_n(X_n; x_{n-1} \leftarrow x_{n-1} x_{n})
                + B_n(X_n; x_{n-1} \leftarrow x_{n-1} x_{n})
\end{eqnarray*}
Substituting from (\ref{first-part-gf}) 
and Lemma \ref{lemma:Bn} gives the result.
\qed

In order to make use of the recurrence of Proposition \ref{thm:alhc-rec}
to prove the generating function (\ref{antiq}) for
anti-lecture hall compositions,
let $A_n(q,s) \triangleq A_n(q,q,q,\ldots,q,s)$.
Then the recurrence of Proposition \ref{thm:alhc-rec} becomes
\begin{eqnarray}
A_n(q,s) = \frac{A_{n-1}(q,q)}{1-s} - A_{n-1}(q,qs)
\frac{s(1-s^{n-1}q^{{n \choose 2}})}{(1-s)(1-s^nq^{{n \choose 2}})},
\label{Aqs}
\end{eqnarray}
with initial condition $A_0(q,s)=1$.
If we were to proceed as with two-rowed plane partitions, we would
(i) ``guess'' the form of $A_n(q,s)$, (ii) prove by induction
that it satisfies (\ref{Aqs}), and then (iii) show that setting
$s=q$ gives (\ref{antiq}).  This would be the easiest proof and it
would give a refinement of the anti-lecture hall generating function,
enumerating solutions according to both the weight and the size of the last
part:
\[
\sum_{\la \in S_{\A_n}} q^{|\la|}s^{\la_n} = A_n(q,qs).
\]
Since we have {\em not} succeeded in guessing $A_n(q,s)$, we follow
a different approach. 
Iterating the recurrence of (\ref{Aqs}) gives:
\begin{eqnarray}
A_n(q,s) =\sum_{i=0}^{n-1} (-1)^i A_{n-1-i}(q,q) s^i q^{{i \choose 2}}
\frac{1-s^{n-i} q^{{n \choose 2}-{i \choose 2}}}
{(s;q)_{i+1} (1-s^n q^{{n \choose 2}})}.
\label{Aqs-solut}
\end{eqnarray}

\noindent
Now, setting $s=q$ gives a recurrence independent of $s$:
\begin{equation}
A_n(q,q)=\sum_{i=0}^{n-1} (-1)^i A_{n-1-i}(q,q)
\frac{q^{i+1\choose 2}-q^{n+1\choose 2}}{(q;q)_{i+1}(1-q^{n+1\choose 2})}.
\label{antiq-rec}
\end{equation}
We show by induction that
the solution to  {\rm (\ref{antiq-rec})} is
\[
A_n(q,q) = \frac{(-q)_n}{(q^2)_n}.
\]
Assume inductively that $A_{n-1-i}=(-q)_{n-1-i}/(q^2)_{n-1-i}$.
Then we need to prove that
\[
\frac{B_n(q)-q^{n+1\choose 2}C_n(q)}{1-q^{n+1\choose 2}}=
\frac{(-q)_n}{(q^2)_n};
\]
with
\[
C_n(q)=\sum_{i=0}^{n-1} (-1)^i
\frac{(-q)_{n-1-i}}{(q^2)_{n-1-i}(q)_{i+1}}
\]
and
\[
B_n(q)=\sum_{i=0}^{n-1} (-1)^i q^{i(i+1)/2}
\frac{(-q)_{n-1-i}}{(q^2)_{n-1-i}(q)_{i+1}}
\]

We will prove that
\[
B_{2n+1}(q)=\frac{(-q)_{2n+1}}{(q^2)_{2n+1}}\ \ \
C_{2n+1}(q)=\frac{(-q)_{2n+1}}{(q^2)_{2n+1}}
\]
\[
B_{2n}(q)=\frac{(-q)_{2n}}{(q^2)_{2n}}-\frac{q^{2n+1\choose
2}}{(q^2)_{2n}} ;\ \ \ \ \
C_{2n}(q)=\frac{(-q)_{2n}}{(q^2)_{2n}}-\frac{1}{(q^2)_{2n}}
\]

Therefore, we need to prove the following identities for $C_n$~:
\begin{equation}
\sum_{i=0}^{2n} (-1)^i \frac{(-q)_{2n-i}}{(q^2)_{2n-i}(q)_{i+1}}=
\frac{(-q)_{2n+1}}{(q^2)_{2n+1}}. \label{odd1}
\end{equation}
\begin{equation}
\sum_{i=0}^{2n-1} (-1)^i
\frac{(-q)_{2n-1-i}}{(q^2)_{2n-1-i}(q)_{i+1}}=
\frac{(-q)_{2n}}{(q^2)_{2n}}-\frac{1}{(q^2)_{2n}}. \label{even1}
\end{equation}

A few $q$-series manipulations show that the two previous equations
are equivalent to:
\begin{equation}
\sum_{j=0}^{n} (-1)^{j}(-1;q)_{j}\left[\begin{array}{c} n\\
j\end{array}\right]_q=(-1)^n
\end{equation}

Recalling that
$$\left[\begin{array}{c} n\\
j\end{array}\right]_q = \frac{(q^{-n})_j(-1)^jq^{nj -
j(j-1)/2}}{(q)_j},
$$
we see that the identity follows from the case $a = -1, c \to
\infty$ of $q$-Chu Vandermonde summation (1.5.2 in \cite{gasper}),
\begin{equation} \label{qchu1}
\sum_{j=0}^n \frac{(a)_j(q^{-n})_j(cq^n/a)^j}{(c)_j(q)_j} =
\frac{(c/a)_n}{(c)_n}.
\end{equation}

Now we need
\begin{equation}
\sum_{i=0}^{2n} (-1)^i q^{i+1\choose 2}
\frac{(-q)_{2n-i}}{(q^2)_{2n-i}(q)_{i+1}}=
\frac{(-q)_{2n+1}}{(q^2)_{2n+1}}. \label{odd2}
\end{equation}
\begin{equation}
\sum_{i=0}^{2n-1} (-1)^i q^{i(i+1)/2}
\frac{(-q)_{2n-1-i}}{(q^2)_{2n-1-i}(q)_{i+1}}=
\frac{(-q)_{2n}}{(q^2)_{2n}}+\frac{q^{2n+1\choose 2}}{(q^2)_{2n}}.
\label{even2} \end{equation}

The same $q$-series manipulations show that the two previous
equations are equivalent to:
\begin{equation}
\sum_{j=0}^{n} (-1)^{j}(-1;q)_{j}\left[\begin{array}{c} n\\
j\end{array}\right]_q q^{{n-j}\choose 2}=(-1)^nq^{n\choose 2}
\end{equation}

\noindent
This follows in a similar way from the ``other" $q$-Chu Vandermonde
summation (1.5.3 in \cite{gasper}),
\begin{equation} \label{qchu2}
\sum_{j=0}^n \frac{(a)_j(q^{-n})_jq^j}{(c)_j(q)_j} =
\frac{a^n(c/a)_n}{(c)_n},
\end{equation}

\noindent
under the substitutions $a= -1, c = 0$.
\qed

%\noindent
%{\em Observation.}
%If we further expand the recurrence (\ref{Aqs-solut}),
%we get a sum over all compositions of $n$:
%let $C(n,t)$ denote the set of compositions of $n$
%into $t$ positive parts.  Then
%\[
%A_n(q,s)= \sum_{t=0}^n (-1)^{n-i} \sum_{a \in C(n,t)}\prod_{j=1}^t
%T_{a_j,n-(a_1+ \cdots a_{j-1})}(q,s),
%\]
%where
%\[
%T_{i,n}(q,s)=  s^i q^{{i \choose 2}}
%\frac{1-s^{n-i} q^{{n \choose 2}-{i \choose 2}}}
%{(s;q)_{i+1} (1-s^n q^{{n \choose 2}})}.
%\]
%This is what we would expect from an inclusion-exclusion
%approach to counting.  We can interpret the sum (\ref{})
%combinatorially as ....

\pagebreak
\section{Lecture Hall Partitions}

In \cite{BME1}, Bousquet-M\'elou and Eriksson studied the
set of integer sequences $\la=(\la_1, \ldots, \la_n)$
satisfying the constraints
\[
\mL_n \ \ \ = \ \ \ \left [
\frac{\la_1}{n}\ge\frac{\la_2}{n-1}\ge \ldots \ge \frac{\la_n}{1}\ge 0
\right ].
\]
They referred to these as
{\em lecture hall partitions}
and showed that the generating function
is
\begin{equation}
L_n(q) \triangleq
\sum_{\lambda\in S_{\mL_n}} q^{|\lambda|} =
\prod_{i=1}^{n}\frac{1}{1-q^{2i-1}}.
\label{lhpq}
\end{equation}
In \cite{PA1}, Andrews showed how to use partition analysis to derive
a recurrence for the full generating function of $\mL_n$.  However,
substantial new ideas, outside of partition analysis, were required to 
move from this to the solution (\ref{lhpq}).  

In this section, we show that by
strategic application of Theorem \ref{thm:guidelines},
we can derive a recurrence
for the full generating function of a generalization of
$\mL_n$ that will reduce the
proof of (\ref{lhpq}) to a 
$q$-series calculation (albeit nontrivial).
Our derivation  here via the five guidelines is
both simpler and more elementary than the approach in \cite{CS2} 
(at the expense of a more challenging $q$-series calculation).

%First consider a generalization.
In \cite{CS2}, we defined
{\em truncated lecture hall partitions} to be the integer sequences
satisfying:
\[
\mL_{n,k} \ \ \ = \ \ \ \left [
\frac{\la_1}{n}\ge\frac{\la_2}{n-1}\ge \ldots \ge \frac{\la_k}{n-k+1}\ge 0
\right ].
\]
We showed that if 
\begin{equation}
\bar{\mL}_{n,k} \ \ \ = \ \ \ \left [
\frac{\la_1}{n}\ge\frac{\la_2}{n-1}\ge \ldots \ge \frac{\la_k}{n-k+1}  >0
\right ],
\label{trunc-constraints}
\end{equation}
that is, all parts must be positive, the generating function is
\begin{equation}
\bar{L}_{n,k}(q)=q^{k+1\choose 2}\left[\begin{array}{l}n\\ k\end{array}\right]
_q
\frac{(-q^{n-k+1};q)_{k}}{(q^{2n-k+1};q)_{k}}.
\label{Lnk}
\end{equation}
It can be checked that
setting $k=n$ and dividing by $q^{{n+1 \choose 2}}$ gives
(\ref{lhpq}).  

%We show how to use Theorem \ref{thm:guidelines} to
%get a recurrence for $\bar{L}_{n,k}(x_1, \ldots,x_k)$ in terms of 
%$\bar{L}_{n,k-1}(x_1, \ldots,x_{k-1})$.

%Let $\bar{L}_{n,k}(x_1,\ldots ,x_k)$ be the gf of LHP in $L_n$ with
%$k$ positive parts. $x_i$ counts the weight of the $ith$ part.

\begin{proposition}
The generating function for truncated lecture hall partitions
{\rm (\ref{trunc-constraints})} satisfies
\begin{eqnarray*}
%\bar{L}_{n,k}(x_1,\ldots
%,x_k)&=&\frac{x_k}{1-x_k}\bar{L}_{n,k-1}(x_1,\ldots ,x_{k-1})-
%\frac{1}{1-x_k}\bar{L}_{n,k-1}(x_1,\ldots ,x_{k-2},x_{k-1}x_k)\\&&-
%\frac{z_{n,k}}{1-z_{n,k}}\bar{L}_{n,k-1}(x_1,\ldots
%,x_{k-2},x_{k-1}x_k)
\bar{L}_{n,k}(x_1,\ldots ,x_k)&=&
\frac{x_k \bar{L}_{n,k-1}(x_1,\ldots ,x_{k-1})}
{1-x_k}-
\frac{
\bar{L}_{n,k-1}(x_1,\ldots ,x_{k-2},x_{k-1}x_k)}
{1-x_k}
\\&&
-
\frac
{z_{n,k}
\bar{L}_{n,k-1}(x_1,\ldots ,x_{k-2},x_{k-1}x_k)}
{1-z_{n,k}}.
 \end{eqnarray*}
with
$z_{n,k}=x_1^nx_2^{n-1}\ldots x_k^{n-k+1}$.
\label{prop:lhp}
\end{proposition}
\noindent{\bf Proof.} 
Note that $\la_{k-1} > \la_k$ is implied by $\bar{\mL}_{n,k}$,
so by guideline 4, 
$\bar{\mL}_{n,k} = \bar{\mL}_{n,k} \cup \{[\la_{k-1} > \la_k] \}$.
Now apply guideline 5 with $c=[\la_{k-1}  \geq  \frac{n-k+2}{n-k+1} \la_{k}]$
to get $\bar{\mL}_{n,k}= \D-\E$:
% $\bar{\mL}_{n,k} = \D-\E$:
\begin{equation}
{\small
\bar{\mL}_{n,k} =  
\left [ \begin{array}{rcl}
\la_1 & \geq  &\frac{n}{n-1} \la_2 \\
\vspace{-.05in}\\
\la_2  &\geq  &\frac{n-1}{n-2} \la_3 \\
&\vdots &\\
\la_{k-3} & \geq  &\frac{n-k+4}{n-k+3} \la_{k-2} \\
\vspace{-.05in}\\
\la_{k-2} & \geq  &\frac{n-k+3}{n-k+2} \la_{k-1} \\
\vspace{-.05in}\\
\la_{k-1}  &\geq  &\frac{n-k+2}{n-k+1} \la_{k} \\
\vspace{-.05in}\\
\la_{k-1}  &> & \la_k\\
\vspace{-.05in}\\
\la_k  &>&0
\end{array} \right ]
= 
\left [ \begin{array}{rcl}
\la_1 & \geq  &\frac{n}{n-1} \la_2 \\
\vspace{-.05in}\\
\la_2  &\geq  &\frac{n-1}{n-2} \la_3 \\
&\vdots &\\
\la_{k-3} & \geq  &\frac{n-k+4}{n-k+3} \la_{k-2} \\
\vspace{-.05in}\\
\la_{k-2} & \geq  &\frac{n-k+3}{n-k+2} \la_{k-1} \\
\vspace{-.05in}\\
\\
\vspace{-.05in}\\
\la_{k-1}  &> & \la_k\\
\vspace{-.05in}\\
\la_k  &>&0
\end{array} \right ]
 -  
\left [ \begin{array}{rcl}
\la_1 & \geq  &\frac{n}{n-1} \la_2 \\
\vspace{-.05in}\\
\la_2  &\geq  &\frac{n-1}{n-2} \la_3 \\
&\vdots &\\
\la_{k-3} & \geq  &\frac{n-k+4}{n-k+3} \la_{k-2} \\
\vspace{-.05in}\\
\la_{k-2} & \geq  &\frac{n-k+3}{n-k+2} \la_{k-1} \\
\vspace{-.05in}\\
\la_{k}  &>  &\frac{n-k+1}{n-k+2} \la_{k-1} \\
\vspace{-.05in}\\
\la_{k-1}  &> & \la_k\\
\vspace{-.05in}\\
\la_k  &>&0
\end{array} \right ]
\label{lhpL}
}
\end{equation}
The first system on the right, $\D$, implies the constraint
$\la_{k-1}>0$, so it can be added.  Now apply guideline 5 to
$\D$ using $c=[\la_{k-1}  >  \la_k]$ to get:
% $\D = \F-\G$.
\begin{equation}
{\small
\D = 
\left [ \begin{array}{rcl}
\la_1 & \geq  &\frac{n}{n-1} \la_2 \\
\vspace{-.05in}\\
\la_2  &\geq  &\frac{n-1}{n-2} \la_3 \\
&\vdots &\\
\la_{k-3} & \geq  &\frac{n-k+4}{n-k+3} \la_{k-2} \\
\vspace{-.05in}\\
\la_{k-2} & \geq  &\frac{n-k+3}{n-k+2} \la_{k-1} \\
\vspace{-.05in}\\
\la_{k-1}  &> & \la_k\\
\vspace{-.05in}\\
\la_{k-1}  &>&0\\
\vspace{-.05in}\\
\la_k  &>&0
\end{array} \right ]
=  
\left [ \begin{array}{rcl}
\la_1 & \geq  &\frac{n}{n-1} \la_2 \\
\vspace{-.05in}\\
\la_2  &\geq  &\frac{n-1}{n-2} \la_3 \\
&\vdots &\\
\la_{k-3} & \geq  &\frac{n-k+4}{n-k+3} \la_{k-2} \\
\vspace{-.05in}\\
\la_{k-2} & \geq  &\frac{n-k+3}{n-k+2} \la_{k-1} \\
\vspace{-.05in}\\
\\
\vspace{-.05in}\\
\la_{k-1}  &>&0\\
\vspace{-.05in}\\
\la_k  &>&0
\end{array} \right ]
 -  
\left [ \begin{array}{rcl}
\la_1 & \geq  &\frac{n}{n-1} \la_2 \\
\vspace{-.05in}\\
\la_2  &\geq  &\frac{n-1}{n-2} \la_3 \\
&\vdots &\\
\la_{k-3} & \geq  &\frac{n-k+4}{n-k+3} \la_{k-2} \\
\vspace{-.05in}\\
\la_{k-2} & \geq  &\frac{n-k+3}{n-k+2} \la_{k-1} \\
\vspace{-.05in}\\
\la_{k}  &\geq & \la_{k-1}\\
\vspace{-.05in}\\
\la_{k-1}  &>&0\\
\vspace{-.05in}\\
\la_k  &>&0
\end{array} \right ].
\label{lhpD}
}
\end{equation}
The first system on the right of (\ref{lhpD}) is just 
$\bar{L}_{n,k-1} \cup \{[\la_{k} > 0]\}$.
The second system on the right  becomes
$\bar{L}_{n,k-1} \cup \{[\la_{k} \geq 0]\}$ after the substitution
$\la_{k} \leftarrow \la_{k}+  \la_{k-1}$. So, by Theorem \ref{thm:guidelines}
and
summarizing so far, we have
\begin{equation}
\bar{L}_{n,k}(x_1, \ldots, x_k)=
\frac{x_k\bar{L}_{n,k-1}(x_1, \ldots, x_{k-1})}{(1-x_k)} -
\frac{\bar{L}_{n,k-1}(x_1, \ldots, x_{k-1}x_k)}{(1-x_k)} -
E(x_1, \ldots, x_k),
\label{Lnkgf}
\end{equation}
where $E(x_1, \ldots, x_k)$ is the generating function
for the last constraint system,
$\E$, in (\ref{lhpL}).
Apply $\la_{k-1} \leftarrow \la_{k-1}+\la_{k}$ to  $\E$ followed by
 $\la_{k} \leftarrow \la_{k}+(n-k+1)\la_{k-1}$ as illustrated
below
$\E \rightarrow \E' \rightarrow \F$:
\begin{equation}
\label{lhpE}
\end{equation}
\begin{equation*}
{\small
%\bar{\E}_{n,k} =
\left [ \begin{array}{rcl}
\la_1 & \geq  &\frac{n}{n-1} \la_2 \\
\vspace{-.05in}\\
\la_2  &\geq  &\frac{n-1}{n-2} \la_3 \\
&\vdots &\\
\la_{k-3} & \geq  &\frac{n-k+4}{n-k+3} \la_{k-2} \\
\vspace{-.05in}\\
\la_{k-2} & \geq  &\frac{n-k+3}{n-k+2} \la_{k-1} \\
\vspace{-.05in}\\
\la_{k}  &>  &\frac{n-k+1}{n-k+2} \la_{k-1} \\
\vspace{-.05in}\\
\la_{k-1}  &> & \la_k\\
\vspace{-.05in}\\
\la_k  &>&0
\end{array} \right ]
\rightarrow
\left [ \begin{array}{rcl}
\la_1 & \geq  &\frac{n}{n-1} \la_2 \\
\vspace{-.05in}\\
\la_2  &\geq  &\frac{n-1}{n-2} \la_3 \\
&\vdots &\\
\la_{k-3} & \geq  &\frac{n-k+4}{n-k+3} \la_{k-2} \\
\vspace{-.05in}\\
\la_{k-2} & \geq  &\frac{n-k+3}{n-k+2} (\la_{k-1}+\la_{k}) \\
\vspace{-.05in}\\
\la_{k}  &>  &({n-k+1}) \la_{k-1} \\
\vspace{-.05in}\\
\la_{k-1}  &> & 0\\
\vspace{-.05in}\\
\la_k  &>&0
\end{array} \right ]
\rightarrow
\left [ \begin{array}{rcl}
\la_1 & \geq  &\frac{n}{n-1} \la_2 \\
\vspace{-.05in}\\
\la_2  &\geq  &\frac{n-1}{n-2} \la_3 \\
&\vdots &\\
\la_{k-3} & \geq  &\frac{n-k+4}{n-k+3} \la_{k-2} \\
\vspace{-.05in}\\
\la_{k-2} & \geq  &\frac{n-k+3}{n-k+2} \la_{k} + (n-k+3)\la_{k-1} \\
\vspace{-.05in}\\
\la_{k}  &>  &0 \\
\vspace{-.05in}\\
\la_{k-1}  &> & 0\\
\vspace{-.05in}\\
\\
\end{array} \right ].
}
\end{equation*}
%If $\E$ is the first constraint system above, $\E'$ the second,
%and $\F$ the last,
By guideline  3,
\[
E(x_1, \ldots, x_k)=
E'(x_1, \ldots, x_{k-1},x_{k-1}x_k),
\]
\[
E'(x_1, \ldots x_k)= 
F(x_1, \ldots, x_{k-2}, x_{k-1}x_k^{n-k+1},x_k),
\]
so
\begin{equation}
E(x_1, \ldots, x_k)=
F(x_1, \ldots, x_{k-2}, x_{k-1}^{n-k+2}x_k^{n-k+1},x_{k-1}x_k).
\label{Egf}
\end{equation}
Finally, starting from $\F$, the last set of constraints in
(\ref{lhpE}),  perform the following
sequence of substitutions
\[
\la_i \leftarrow \la_i + (n-i+1)\la_{k-1}; \ \ \ i = k-2, \ldots 1,
\]
as illustrated below:
%\begin{equation*}
{\small
\begin{eqnarray}
{\F} & \rightarrow &
\left [ \begin{array}{rcl}
\la_1 & \geq  &\frac{n}{n-1} \la_2 \\
\vspace{-.05in}\\
\la_2  &\geq  &\frac{n-1}{n-2} \la_3 \\
&\vdots &\\
\la_{k-3} & \geq  &\frac{n-k+4}{n-k+3} \la_{k-2} + (n-k+4)\la_{k-1} \\
\vspace{-.05in}\\
\la_{k-2} & \geq  &\frac{n-k+3}{n-k+2} \la_{k}  \\
\vspace{-.05in}\\
\la_{k}  &>  &0 \\
\vspace{-.05in}\\
\la_{k-1}  &> & 0\\
\end{array} \right ]
\rightarrow \ldots \rightarrow
\left [ \begin{array}{rcl}
\la_1 & \geq  &\frac{n}{n-1} \la_2 \\
\vspace{-.05in}\\
\la_2  &\geq  &\frac{n-1}{n-2} \la_3 + (n-1) \la_{k-1} \\
&\vdots &\\
\la_{k-3} & \geq  &\frac{n-k+4}{n-k+3} \la_{k-2}  \\
\vspace{-.05in}\\
\la_{k-2} & \geq  &\frac{n-k+3}{n-k+2} \la_{k}  \\
\vspace{-.05in}\\
\la_{k}  &>  &0 \\
\vspace{-.05in}\\
\la_{k-1}  &> & 0\\
\end{array} \right ] \nonumber\\
&&
\longrightarrow 
\left [ \begin{array}{rcl}
\la_1 & \geq  &\frac{n}{n-1} \la_2 + n \la_{k-1} \\
\vspace{-.05in}\\
\la_2  &\geq  &\frac{n-1}{n-2} \la_3   \\
&\vdots &\\
\la_{k-3} & \geq  &\frac{n-k+4}{n-k+3} \la_{k-2}  \\
\vspace{-.05in}\\
\la_{k-2} & \geq  &\frac{n-k+3}{n-k+2} \la_{k}  \\
\vspace{-.05in}\\
\la_{k}  &>  &0 \\
\vspace{-.05in}\\
\la_{k-1}  &> & 0\\
\end{array} \right ]
\ \ \ \longrightarrow \ \ \ 
\left [ \begin{array}{rcl}
\la_1 & \geq  &\frac{n}{n-1} \la_2   \\
\vspace{-.05in}\\
\la_2  &\geq  &\frac{n-1}{n-2} \la_3   \\
&\vdots &\\
\la_{k-3} & \geq  &\frac{n-k+4}{n-k+3} \la_{k-2}  \\
\vspace{-.05in}\\
\la_{k-2} & \geq  &\frac{n-k+3}{n-k+2} \la_{k}  \\
\vspace{-.05in}\\
\la_{k}  &>  &0 \\
\vspace{-.05in}\\
\la_{k-1}  &> & 0\\
\end{array} \right ] \ \ \ = \ \ \ \G.
\label{lhpF}
\end{eqnarray}
}
The resulting system of constraints, $\G$, in (\ref{lhpF})
can be viewed as
$\mL_{n,k-1}$, where $\la_{k-1}$ has been replaced by $\la_{k}$,
together with the constraint
$[\la_{k-1} > 0]$. Thus,
\begin{equation}
G(x_1, \ldots, x_k) =
\frac{x_{k-1}L_{n,k-1}(x_1, \ldots, x_{k-2},x_k)}{(1-x_{k-1})}.
\label{Ggf}
\end{equation}
By guideline 3, the generating function for $\F$
is obtained from $G$  
by the sequence of substitutions
\[
x_{k-1} \leftarrow x_{k-1} x_{i}^{n-i+1}; \ \ \ i = 1 \ldots k-2,
\]
giving
\begin{equation}
F(x_1, \ldots, x_k) =
G(x_1, \ldots, x_{k-2}, x_1^nx_2^{n-1} \cdots x_{k-2}^{n-k+3}x_{k-1}, x_k).
\label{Fgf}
\end{equation}
Returning to $E$ in (\ref{Egf}) and using (\ref{Ggf}) and (\ref{Fgf}),
\begin{eqnarray}
E(x_1, \ldots, x_k) &=&
F(x_1, \ldots, x_{k-2}, x_{k-1}^{n-k+2}x_k^{n-k+1},x_{k-1}x_k)
\nonumber\\
& = & 
G(x_1,\ldots,x_{k-2},x_1^nx_2^{n-1} \cdots x_{k-2}^{n-k+3} x_{k-1}^{n-k+2}
x_k^{n-k+1}, x_{k-1}x_k)
\nonumber \\
& = & 
\frac{x_1^nx_2^{n-1} \cdots   x_k^{n-k+1}
L_{n,k-1}(x_1,\ldots,x_{k-2}, x_{k-1}x_k)}
{1-x_1^nx_2^{n-1} \cdots x_k^{n-k+1}}.
\label{Efull}
\end{eqnarray}
Combining (\ref{Efull})  with (\ref{Lnkgf}) gives the result.
\qed

Let 
$\bar{L}_{n,k}(q,s) = \bar{L}_{n,k}(q,q, \ldots, q,s).$
Setting $x_k=s$ and $x_i=q$ for $i<k$ in Proposition \ref{prop:lhp}
gives 
\begin{equation}
\bar{L}_{n,k}(q,s)=\frac{s}{1-s}\bar{L}_{n,k-1}(q,q)-
\bar{L}_{n,k-1}(q,sq)
\left (
\frac{1}{1-s}+\frac{z_{n,k}}{1-z_{n,k}}
\right ),
\label{Lnkqsrec}
\end{equation}
where $z_{n,k}=s^{n-k+1}q^{{n+1\choose 2}-{n-k+2\choose 2}}$.
One would hope to prove (\ref{Lnk}) now by finding a closed form for
$\bar{L}_{n,k}(q,s)$, proving that it satisfies the recurrence
(\ref{Lnkqsrec}) and then setting $s=q$ to get
(\ref{Lnk}). Since
we were unable to
guess $\bar{L}_{n,k}(q,s)$, we proceed as for
anti-lecture hall compositions  to iterate the recurrence (\ref{Lnkqsrec})
and get
$$
\bar{L}_{n,k}(q,s)=\sum_{j\ge 1}(-1)^{j-1}\frac{sq^{j-1}}{(s;q)_j}\cdot
\frac{1-s^{n-k+j}q^{(n-k+j)(j-2)+{n+1\choose 2}-{n-k+j\choose 2}}}{1-s^{n-k+1}
q^{{n+1\choose 2}-{n-k+2\choose 2}}}\cdot \bar{L}_{n,k-j}(q,q).
$$
Now setting $s=q$ we need only a single argument:
$$
\bar{L}_{n,k}(q)=\sum_{j\ge 1}(-1)^{j-1}\frac{q^{j}}{(q;q)_j}\cdot
\frac{1-q^{k(n-k+j)+{k-j+1\choose 2}}}{1-q^{{n+1\choose 2}-{n-k+1\choose 2}}}
\cdot \bar{L}_{n,k-j}(q).
$$
It remains to prove that this recurrence is satisfied by (\ref{Lnk}).
We defer the details until a later report;  our main point was to show
that strategic application of the guidelines reduce the truncated
lecture hall theorem to a $q$-series computation.

\section{The Five Guidelines Suffice}
\label{section:suffice}

Let $\C$ be the set of inequalities
\begin{equation}
c_{i,0} +  c_{i,1}\lambda_1 +
c_{i,2}\lambda_2 + \ldots
+ c_{i,n}\lambda_n \geq  0, \ \ \
 1 \leq i \leq r.
\label{suffice-constraints}
\end{equation}
and let $S_{\C}$ be the set of
of nonnegative integer sequences
satisfying all constraints in  $\C$.
In this section we show that the five guidelines of Theorem 
\ref{thm:guidelines} are powerful enough to find the generating
function of  $S_{\C}$
for any integers $c_{i,j}$.
We will assume that all constraints are {\em homogeneous},
i.e., that $c_{i,0}=0$.
Otherwise, 
introduce a new variable $\la_0$ and let $\C'$ be the same as $\C$,
except that for every $i$, the $i$th constraint is now: 
\[
c_{i,0}\la_0 +  c_{i,1}\lambda_1 +
c_{i,2}\lambda_2 + \ldots
+ c_{i,n}\lambda_n \geq  0.
\]
Then 
$F_{\C}(x_1, \ldots x_n)$  is the coefficient of 
$x_0$ in $F_{\C'}(x_0,x_1, \ldots x_n)$. 
We also generalize the claim a bit to allow any of the constraints of $\C$ to be
equalities.
\begin{theorem}
The five guidelines of Theorem \ref{thm:guidelines}
are sufficient to find the full generating function
for any homogeneous system of linear inequalities and equalities.
\label{thm:suffice}
\end{theorem}
\noindent
{\bf Proof.}
Let $\C$ be a homogeneous system of linear inequalities and equalities
 with variables
$\la_1, \ldots, \la_n$.  Since we require nonnegative integer solutions,
we can assume that for each variable $\la_i$, $\C$ contains a constraint
$b_i$ of the form $[\la_i \geq 0]$ or $[\la_i = 0]$. Call these constraints
$b_i$ {\em basic}.
Write $\C$ in the form
\[
\C = [c_1,c_2, \ldots, c_r;b_1,b_2, \ldots, b_n],
\]
where $c_1,c_2, \ldots, c_r$ is an ordered list of the non-basic
constraints in $\C$.
If $r=0$, all constraints are basic and the generating function follows
from guidelines 1 and 2 (and $F_{[\la_i=0]}(x_i)=1$).

%Otherwise, let
%\begin{quote}
%\vspace{-.3in}
%\item
%$M$ be the largest positive coefficient of $c_1$ (0, if none);
%\item
%$e_{max}$ be the number of occurrences of $M$ among the coefficients of $c_1$;
%\item
%$m$ be the smallest negative coefficient of $c_1$ (0, if none);
%\item
%$e_{min}$ be the number of occurrences of $m$ among the coefficients of $c_1$.
%\end{quote}
%When $r>0$ we show that we can use the guidelines to reduce the
%computation of the generating function of $\C$ to the computation of
%the generating function of one or more systems $\C'$ in which 
%at least one of the statistics $\{r,M,e_{max},|m|,e_{min}\}$ has been reduced.

\noindent
Otherwise, define:

$M$:  the largest positive coefficient of $c_1$ (0, if none);

$e_{max}$:  the number of occurrences of $M$ among the coefficients of $c_1$;

$m$:  the smallest negative coefficient of $c_1$ (0, if none);

$e_{min}$:  the number of occurrences of $m$ among the coefficients of $c_1$.

\noindent
When $r>0$ we show that we can use the guidelines to reduce the
computation of the generating function of $\C$ to the computation of
the generating function of one or more systems $\C'$ in which
at least one of the statistics $\{r,M,e_{max},|m|,e_{min}\}$ has been reduced.

If $m=0$, all coefficients of $c_1$ are nonnegative, so $c_1$ is
redundant and can be deleted.
Otherwise, if $M=0$, all coefficients of $c_1$ are nonpositive and so we get
an equivalient system replacing $\la_j$ by 0 in $c_1, \ldots, c_r$ and
setting $b_j=[\la_j=0]$.
In so doing we have decreased $|m|$ or $e_{min}$.

Otherwise, $m<0$ and $M>0$;  we do a version of {\em Elliott reduction}
\cite{elliott}.  Let $i$ and $j$ be such that
$m$ is the coefficient of $\la_i$ in $c_1$ and
$M$ is the coefficient of $\la_j$.
We would like to use guideline 3 and reduce to a system with smaller
$M$ or $e_{max}$ or $|m|$ or $e_{min}$.
First use guideline 4 with $c=[\la_i \geq \la_j]$:
\[
F_{\C}(X_n) = 
F_{\C \cup [\la_i \geq \la_j]}(X_n) + 
F_{\C \cup [\la_j > \la_i]}(X_n). 
\]
For the first term, $F_{\C \cup [\la_i \geq \la_j]}(X_n)$, do the 
substitution
$\la_i \leftarrow \la_i +\la_j $ into constraints 
$c_1,c_2, \ldots, c_r$ in $\C$.  This decreases the coefficient of $\la_j$,
thereby decreasing $M$ or $e_{max}$.  By guideline 3, the substitution
$x_j \leftarrow x_jx_i$ in the generating function of the resulting
constraint system gives $F_{\C \cup [\la_i \geq \la_j]}(X_n)$.

For the second term, $F_{\C \cup [\la_j > \la_i]}(X_n)$, if
$b_j = [\la_j=0]$, there are no solutions and the generating function is 0.
Otherwise, $b_j = [\la_j \geq 0]$.
Substitute $\la_j \leftarrow \la_j +\la_i $ into constraints
$c_1,c_2, \ldots, c_r$ in $\C$ to get $\C'$.   This increases the
coefficient of $\la_i$, thereby decreasing $|m|$ or $e_{min}$.
Substituting $\la_j \leftarrow \la_j +\la_i $ into $[\la_j > \la_i]$
gives $[\la_j >0]$.  By guideline 3,
\[
F_{\C \cup [\la_j > \la_i]}(X_n)=
F_{\C' \cup [\la_j > 0]}(X_n; x_i \leftarrow x_ix_j).
\]
However, we disallow strict inequalities.  So, use guideline 5 with
$c=[\la_j > 0]$ and observe that $b_j=[\la_j \geq 0] \in \C'$.
Let $\C''$ denote  $\C'$ with $b_j \leftarrow [\la_j=0]$.  Then
\[
F_{\C' \cup [\la_j > 0]}(X_n) =
F_{\C'}(X_n) -
F_{\C' \cup [\la_j \leq 0]}(X_n) =
F_{\C'}(X_n) -
F_{\C''}(X_n). 
\]
\qed

\noindent
(Note that we have optimized the proof for simplicity at the expense of 
algorithmic efficiency.)
  %Nevertheless, note in the last step that
%\[
%F_{\C''}(X_n)= 
%F_{\C'}(X_n; x_j \leftarrow 0). 
%\]
%An equivalent observation  appearing in \cite{} is used to speed up
%computation in the Omega package.

It follows from Theorem \ref{thm:suffice} and its proof that the full generating
function of (\ref{suffice-constraints}) can be built up from the functions 
$1/(1-x_i)$ by a finite number of additions, subtractions, and
substitutions.  We get then as a corollary the following well-known
result:
The full generating function for the nonnegative integer solutions to
any system of linear inequalities in $n$ variables with integer
coefficients has the form
\[
\frac{p(x_1, \ldots, x_n)}
{(1-\alpha_1)(1-\alpha_2) \cdots (1-\alpha_t)},
\]
where $t \geq 0$, $p$ is a polynomial in $x_1, \ldots, x_n$
and each $\alpha_i$ is a monomial in $x_1, \ldots, x_n$.

%\begin{corollary}
%The full generating function for the nonnegative integer solutions to
%any system of linear inequalities in $n$ variables with integer
%coefficients has the form
%\[
%\frac{p(x_1, \ldots, x_n)}
%{(1-\alpha_1)(1-\alpha_2) \cdots (1-\alpha_t)},
%\]
%where $t \geq 0$, $p$ is a polynomial in $x_1, \ldots, x_n$
%and each $\alpha_i$ is a monomial in $x_1, \ldots, x_n$.
%\end{corollary}

\section{Relationship to MacMahon's Partition Analysis}

We give a brief introduction to partition analysis
in order to highlight the fact that the guidelines of Theorem 
\ref{thm:guidelines} underlie
the work of MacMahon.
Indeed, they were
distilled from partition analysis by a study of MacMahon's work in \cite{Mamo2}
and its application by Andrews, Paule and Riese 
in the series
of papers \cite{PA1,PA2,PA3,PA4,PA5,PA6,PA7,PA8,PA9,PA10,PA11}.

Consider the
set of constraints $\C = \{c_1, \ldots c_r\}$ where
\[
c_i = [a_{i,1}\la_1 + \cdots +  a_{i,n} \la_n \geq 0].
\]
We seek the full generating function for the set
$S_{\C}$ of nonnegative integer sequences
$\la = (\la_1, \ldots, \la_n)$  satisfying $c_i$,
$1 \leq i \leq r$:
\[
F_{\C}(x_1, \ldots, x_n) = 
\sum_{\la \in S_{\C}} x_1^{\la_1} x_2^{\la_2} \cdots x_n^{\la_n}.
\]
The method of {\em partition analysis}, developed by MacMahon
in \cite{Mamo2} is to view the problem as follows. Let
\begin{equation}
P(x_1,x_2, \ldots, x_n,c_1,c_2, \ldots, c_n) =
\prod_{j=1}^{n}
\frac{1}
{1- x_j c_1^{a_{1,j}} c_2^{a_{2,j}} \cdots  c_r^{a_{r,j}}}.
\label{PAproduct}
\end{equation}
Expanding $P$ gives
\begin{eqnarray}
P(x_1,x_2, \ldots, x_n,c_1,c_2, \ldots, c_n) & = &
\sum_{\la_1, \ldots \la_n \geq 0}
\prod_{j=1}^{n}(x_j c_1^{a_{1,j}} c_2^{a_{2,j}} \cdots  c_r^{a_{r,j}})^{\la_j}
\nonumber \\
&= & \sum_{\la_1, \ldots \la_n \geq 0}
\left (
x_1^{\la_1}x_2^{\la_2} \cdots x_n^{\la_n}
\prod_{i=1}^{r}c_i^{a_{i,1}\la_1 + \cdots  a_{i,n}\la_n} 
\right ).
\label{PAsum}
\end{eqnarray}
Observe that  $\la \in S_{\C}$ iff in the term corresponding to
$\la$ in the sum (\ref{PAsum}) every $c_i$, $1 \leq i \leq r$,
has nonnegative exponent. 
Thus 
$F_{\C}(x_1, \ldots, x_n)$ is recovered from
$P(x_1,x_2, \ldots, x_n,c_1,c_2, \ldots, c_n)$
by deleting all terms in which some $c_i$ has a negative
exponent and then setting $c_1=c_2 = \cdots = c_r = 1$.
MacMahon uses the { \em Omega operator}  to express this
process:
\[
F_{\C}(x_1, \ldots, x_n)  =  \ \ 
\mathrel{\mathop{\Omega}\limits_{\ge}}
P(x_1,x_2, \ldots, x_n,c_1,c_2, \ldots, c_n).
\]
The core of partition analysis is a system of
$Omega$-rules designed to be applied strategically to
transform $P(x_1,x_2, \ldots, x_n,c_1,c_2, \ldots, c_n)$
step-by-step into
$F_{\C}(x_1, \ldots, x_n)$.
This view converts the combinatorial problem into
an algebraic one, opening the possibility, for example, of
a partial fraction decomposition of (\ref{PAproduct})
to assist in the transformation from $P$ to $F$.
A list of basic Omega-rules appears in \cite{Mamo2}(pp. 103-106) and \cite{PA1}.

This approach has proven both powerful and systematic in the
computer solution of systems of inequalities.
However, for deriving recurrences for infinite
families, we found that a
return to some of the basic underlying ideas 
simplified the process. 
We note the roots of guidelines 3-5 of Theorem \ref{thm:guidelines}
 in the work of MacMahon \cite{Mamo2}.

Our use of
guideline 3 (which performs limited column operations on the
constraint matrix) is used to much the same effect as the following
Omega-rule:

\[
\mathrel{\mathop{\Omega}\limits_{\ge}}  \
\frac{1}{(1-x_ic)(1-\frac{x_j}{c^a})} =
%\mathrel{\mathop{\Omega}\limits_{\ge}}  \
\frac{1}{(1-x_i)(1-x_j x_i^a)}.
\]

The utility of guideline 4 was recognized by  MacMahon.
He writes
in \cite{Mamo2}, p. 103,
``A very useful principle is that of
adding an inequality which is a\` fortiori true."
It is also used
in a decomposition shown at the beginning of
\cite{Mamo2}, Section 379, p. 131.
Guideline 5  {\em is} one of MacMahon's   Omega-rules, found in
 \cite{Mamo2}, Section 351, p. 104 (slightly transformed):
\[
\mathrel{\mathop{\Omega}\limits_{\ge}}  \
P(c) = P(1)-
\mathrel{\mathop{\Omega}\limits_{\ge}}  \
P(1/c).  
\]

%\section{A Note on Efficient Implementations}
%
%The recurrences  arising from the techniques of this paper
%result in programs which can be implemented for 
%computational/experimental purposes.  In order to make them feasible
%for larger problem sizes, the standard technique of dynamic programming,
%with or without memoizing 
%is preferable to a recursive implementation
%(see for example \cite{CLR}). 
%For example, consider a recurrence of the form:
%\begin{equation}
%F_n(q,s) = g(n,q,s)F_{n-1}(q,q) + h(n,q,s)F_{n-1}(q,qs),
%\label{efficient}
%\end{equation}
%with initial condition $F_1(q,s)= f(q,s)$, which appears to require
%two recursive calls to $F_{n-1}$.  The method of dynamic programming is
%to compute and store $F_i(q,s)$ for $i = 1, \ldots, n$, say in an array
%$A[1..n]$ as follows:
%$F_1(q,s)$, which is given, is stored in $A[1]$.
%Then for $1 < i \leq n$, $F_i(q,s)$ is computed by looking up
%$F_{i-1}(q,s)$ in $A[i-1]$, performing the substitutions,
%multiplications, and additions indicated by 
%(\ref{efficient}), and then storing the result in $A[i]$. 
%In this way, in the computation of $F_n$,
 %the number of ``calls'' to $F$ becomes linear in $n$ rather than
%exponential in $n$.

\section{Concluding Remarks}

The ``five guidelines'' of Theorem 1 provide a unified setting for 
computing the full generating
function for many challenging families of constraints.
 %non
 %(and weight-generating function) for
%many known things, including
%Two-rowed plane partitions, anti-lecture hall compositions,
%truncated lecture hall partitions, truncated anti-lecture hall compositions,
%Cayley compositions, Carlitz compositions, $d$-compositions, and
%up-down compositions,  (not to mention
%%ordinary partitions and compositions, partitions into distinct
%parts, in fact, all of the ones we can handle with the "C-matrix"
%technique, since anything we can do that way we can also do
%using only guidelines 1,2, amd 3) and suggests 
%new ones (e.g. the ones Sunyoung is
%looking at, rotated anti-lecture hall, others)
However, even though they
are guaranteed to be sufficient to find the generating function for
any homogeneous linear system, we are not necessarily guaranteed
to be able to use them to devise a {\em recurrence} for
a parametrized family of constraint sets.  
%We have a few we don't
%know what to do with yet.
%
%Even when we can't solve the recurrences, we get a program.

%Even though the Five Guidelines ``suffice'', by no means
%do we suggest using Theorem \ref{thm:suffice}
%for computation purposes because of the
%uncontrolled branching.
%The Omega package includes "Omega rules" particularly
%tailored to reduce branching and new work in the area is
%pursues this further (cite Xin and the other paper - someone's student?)

In continuing work
we consider the case when
all constraints have the form $\la_i \geq \la_j$ or
$\la_i > \la_j$, forming
a directed graph.  We show how to get a recurrence by strategically
manipulating the diagrams.
Many examples are presented, including
two- and three-rowed plane partitions, plane partitions with
diagonals, plane partition diamonds, and hexagonal plane partitions.

Finally, we note that in \cite{MR2097324}, Xin offers a speed-up to the
Omega package for implementing MacMahon's partition analysis.
Xin's method uses the theory of iterated Laurent series and partial
fraction decompositions.

\vspace{.2in}
\noindent
{\bf Acknowledgement.}
We are grateful to the referee for a careful reading of the
manuscript and detailed
suggestions to improve the presentation.

{\small
\bibliography{inequalities}
\bibliographystyle{plain}
}

\end{document}